\documentclass[12pt,reqno]{amsart}

\usepackage[lite]{amsrefs}

\usepackage{amssymb}
\usepackage[hidelinks]{hyperref}

\usepackage[all,cmtip]{xy}

\usepackage[shortlabels]{enumitem}

\usepackage{geometry}
\usepackage{cases}

\newcommand{\N}{\mathbb{N}}
\newcommand{\Z}{\mathbb{Z}}
\newcommand{\R}{\mathbb{R}}

\newcommand{\pa}[1]{\left(#1\right)}
\newcommand{\pb}[1]{\left[#1\right]}
\newcommand{\pc}[1]{\left\{#1\right\}}
\newcommand{\pd}[1]{\left|#1\right|}
\newcommand{\pe}[1]{\left\lfloor#1\right\rfloor}

\numberwithin{equation}{section}

\theoremstyle{plain} 
\newtheorem{thm}{Theorem}[section]
\newtheorem{cor}[thm]{Corollary}
\newtheorem{lem}[thm]{Lemma}
\newtheorem{prop}[thm]{Proposition}

\theoremstyle{definition}

\theoremstyle{remark}
\newtheorem{rem}[thm]{Remark}

\begin{document}

\title[Integer Solutions of Pell Equation in Bounded Regions]{Integer Solutions of Pell Equation in Bounded Regions}

\author{ONG KUN YI}
\author{Eddie Shahril Bin Ismail}
\address{Department of Mathematical Sciences, Faculty of Science and Technology, Universiti Kebangsaan Malaysia, Bangi, Selangor, Malaysia}


\email{kunyi070101@gmail.com, esbi@ukm.edu.my}



\date{}

\begin{abstract}
    The Pell equation \(x^2 - Dy^2 = 1\) with non-square \(D > 1\) has infinitely many integer solutions, yet most research has centered on the asymptotic behavior of fundamental units as \(D\) varies. By contrast, the exact distribution of solutions for a fixed \(D\) within bounded regions has received little attention. In this paper, we contribute to this direction by giving an explicit enumeration of all solutions to the Pell equation inside the square \(|x| + |y| \leq \lambda\) for any \(\lambda > 0\). We further extend our results to the \textit{shifted Pell equation} \(\pa{x-a}^2 - D\pa{y-b}^2 = 1\) for integers \(a\) and \(b\), obtaining exact counts for sufficiently large \(\lambda\).
\end{abstract}

\maketitle

\section{Introduction}
Let \(D\) be a non-square positive integer. The Pell equation is the Diophantine equation
\begin{align}
    \label{1.1}
    x^2 - Dy^2 = 1
\end{align}
whose solutions are integer pairs \((x, y)\) satisfying \eqref{1.1}. Its history dates back to the ancient Greeks, including Archimedes’ Cattle Problem (see \cites{Dorrie, Vardi, Dickson}) and to the work of Brahmagupta and Bhaskara in India, as well as Fermat and Euler in Europe. It was Lagrange who finally established the fundamental fact that Pell equation \eqref{1.1} possesses infinitely many integer solutions. Proofs of this result can be found in numerous books and articles on number theory, see for example \cites{Cohn, Barbeau, Lenstra}. The positive solution \((x,y) \in \N^2\) with the smallest value of \(x + y\sqrt{D}\) among them is called the \emph{fundamental solution} and is denoted by \((\alpha, \beta)\). In particular, it can be determined by
\[\alpha + \beta \sqrt{D} = \inf\pc{x + y\sqrt{D}: \pa{x,y} \in \Z^2, x^2 - Dy^2 = 1, x + y\sqrt{D} > 1}.\]

The following well-known result describes the set of integer solutions to Equation \eqref{1.1} (see \cites{Weil, Jacobson, Lehmer, Lenstra, Nagell, LeVeque, Adler}). \\

\begin{prop}
    \label{prop:1.1}
    Let \(F_D\) be defined as the set of all integer solutions to \eqref{1.1}, i.e.
    \begin{align}
        \label{1.2}
        F_D := \pc{(x,y) \in \Z^2: x^2 - Dy^2 = 1}.
    \end{align}
    Then
    \begin{align*}
        F_D = \bigcup_{i,j \in \pc{1,-1}} \pc{\pa{iu_n, jv_n} \in \Z^2: n \in \Z_{\geq 0}},
    \end{align*}
    where \(u_n\) and \(v_n\) are integers given by \(u_n + v_n\sqrt{D} = \pa{\alpha + \beta\sqrt{D}}^n\) for all integers \(n \geq 0\). \\
\end{prop}

Define functions \(u:[0,\infty) \to [1,\infty)\) and \(v:[0,\infty) \to [0,\infty)\) by
\begin{align}
    \label{1.3}
    u(x) &:= \dfrac{1}{2}\pc{\pa{\alpha + \beta\sqrt{D}}^x + \pa{\alpha - \beta\sqrt{D}}^x} \\
    \label{1.4}
    v(x) &:= \dfrac{1}{2\sqrt{D}}\pc{\pa{\alpha + \beta\sqrt{D}}^x - \pa{\alpha - \beta\sqrt{D}}^x}
\end{align}
where \((\alpha, \beta)\) is the fundamental solution. We now state an alternative result to Proposition \ref{prop:1.1} describing the set of integer solutions to Equation \eqref{1.1} in terms of the functions \(u\) and \(v\). \\

\begin{prop}[\cite{Burton}]
    \label{prop:1.2}
    Consider the functions \(u\) and \(v\) as defined in \eqref{1.3} and \eqref{1.4}, respectively. Then \(u(n) = u_n\) and \(v(n) = v_n\) for all integers \(n \geq 0\), where \(u_n\) and \(v_n\) are integers as described in Proposition \ref{prop:1.1}. This implies
    \begin{align*}
        F_D = \bigcup_{i,j \in \pc{1,-1}} \pc{(iu(x), jv(x)) \in \Z^2: x \in \Z_{\geq 0}}.
    \end{align*}
\end{prop}
\vspace{10pt}

Much of the literature has focused on the growth of the \textit{fundamental unit} \(\varepsilon_D := \alpha + \beta\sqrt{D}\) associated with the fundamental solution \((\alpha, \beta)\) as the integer \(D\) varies. For example, in 2011 Fouvry and Jouve \cite{Fouvry} proved that the set of parameters \(D \leq x\) for which \(\log \varepsilon_D\) is larger than \(D^{\frac{1}{4}}\) has a cardinality essentially larger than \(x^{\frac{1}{4}} \log^2 x\). More recently, Xi \cite{Xi} established uniform lower bounds for the counting function
\[S^f\pa{x,\alpha} := \left| \pc{\pa{\varepsilon_D, D}: 2 \leq D \leq x, D \neq \square, \varepsilon_D \leq D^{\frac{1}{2}+\alpha}} \right|,\]
proving that for any fixed \(\alpha \in \pb{\dfrac{1}{2}, 1}\), \(S^f\pa{x, \alpha} \gg \sqrt{x}\log^2 x\) as \(x \to +\infty\). These works emphasize the asymptotic distribution of the fundamental unit as \(D\) varies. \\

By contrast, comparatively little attention has been given to the case where \(D\) is fixed and one studies the exact enumeration of integer solutions in bounded regions of the plane. To our knowledge, this problem has not been investigated in the literature.

\section{Main Results}

Consider the functions \(u\) and \(v\) as in \eqref{1.3} and \eqref{1.4}, respectively. Define function \(f: [0,\infty) \to [1,\infty)\) by \(f := u + v\). Since \(u\)  and \(v\) are strictly increasing and \(u(0) + v(0) = 1\), it follows that \(f\) is bijective and strictly increasing. Therefore, we may define \(g: [1,\infty) \to [0,\infty)\) as the inverse function of \(f\), namely \(g := f^{-1}\), which is also bijective and strictly increasing. \\

In this paper, we study the set of integer solutions to Equation \eqref{1.1}, where \(D\) is a fixed non-square positive integer, lying in the square \(|x| + |y| \leq \lambda\) with \(\lambda > 0\). For any set \(G \subseteq \Z^2\) and any real number \(\lambda > 0\), we define
\[Q_G(\lambda) := \pc{(x,y) \in G: |x| + |y| \leq \lambda}.\]

We now address the problem of enumerating the set \(Q_{F_D}(\lambda)\). As it turns out, the following result holds. \\

\begin{thm}
    \label{thm:2.1}
    Consider the set \(F_D\) as defined in \eqref{1.2}. Then for any real number \(\lambda > 0\),
    \[Q_{F_D}(\lambda) = \left\{
    \begin{aligned}
        &\varnothing \ & \ &\text{if} \ \lambda < 1 \\
        &\bigcup_{i,j \in \pc{1,-1}} \pc{(iu(x), jv(x)) \in \Z^2: x \in \Z_{\geq 0}, x \leq g(\lambda)} \ & \ &\text{if} \ \lambda \geq 1,
    \end{aligned}
    \right.
    \]
    where \(g = f^{-1}\), \(f = u + v\). This implies
    \[\left|Q_{F_D}(\lambda)\right| = \left\{
    \begin{aligned}
        &0 \ & \ &\text{if} \ \lambda < 1 \\
        &2 + 4 \pe{g(\lambda)} \ & \ &\text{if} \ \lambda \geq 1.
    \end{aligned}
    \right.
    \]
\end{thm}
\vspace{10pt}

\begin{proof}
The case of \(\lambda < 1\) is trivial, hence assume that \(\lambda \geq 1\). By Proposition \ref{prop:1.2},
\begin{align*}
    Q_{F_D}(\lambda) &= \bigcup_{i,j \in \pc{1,-1}} \pc{(iu(x), jv(x)) \in \Z^2: x \in \Z_{\geq 0}, \pd{iu(x)} + \pd{jv(x)} \leq \lambda} \\
    &= \bigcup_{i,j \in \pc{1,-1}} \pc{(iu(x), jv(x)) \in \Z^2: x \in \Z_{\geq 0}, f(x) \leq \lambda},
\end{align*}
where \(f = u+v\). Since \(g\) is strictly increasing, \(f(x) \leq \lambda\) implies \(x \leq g\pa{\lambda}\) and the formula of \(Q_{F_D}\pa{\lambda}\) follows. \\

For the number of elements, \(\left|Q_{F_D}(\lambda)\right|\), the case of \(\lambda < 1\) is obvious. For \(\lambda \geq 1\),
\begin{align*}
    Q_{F_D}(\lambda) &= \bigcup_{i,j \in \pc{1,-1}} \pc{(iu(x), jv(x)) \in \Z^2: x \in \Z_{\geq 0}, x \leq g(\lambda)} \\
    &= \pc{(1,0), (-1,0)} \cup \bigcup_{i,j \in \pc{1,-1}} \pc{(iu(x), jv(x)) \in \Z^2: x \in \N, x \leq \pe{g(\lambda)}}
\end{align*}
which implies \(\left| Q_{F_D}(\lambda) \right| = 2 + 4\pe{g(\lambda)}\).
\end{proof}
\vspace{10pt}

Based on Theorem \ref{thm:2.1}, it is natural to ask whether there exists an explicit formula for the function \(g(\lambda)\), where \(\lambda \geq 1\). This will be addressed in the following result, Theorem \ref{thm:2.2}. \\

\begin{thm}
    \label{thm:2.2}
    Consider functions \(u\) and \(v\) as defined in \eqref{1.3} and \eqref{1.4}, respectively. Let \(f = u+v\) and \(g = f^{-1}\). Then for any real number \(x\geq 1\),
    \begin{align}
        \label{2.1}
        g(x) = \dfrac{1}{\log\pa{\alpha + \beta \sqrt{D}}}\cosh^{-1}\pa{\dfrac{Dx - \sqrt{Dx^2 - D + 1}}{D-1}},
    \end{align}
    where \(\cosh^{-1}\) is the inverse cosh function. Equivalently,
    \begin{align}
        \label{2.2}
        g(x) = \dfrac{\log\pa{Dx - \sqrt{Dx^2 - D + 1} + \sqrt{D}\sqrt{\pa{D+1}x^2 - D + 1 - 2x\sqrt{Dx^2 - D + 1}}} - \log\pa{D-1}}{\log\pa{\alpha + \beta \sqrt{D}}}.
    \end{align}
\end{thm}
\vspace{10pt}

\begin{proof}
We begin with function \(f: [0,\infty) \to [1,\infty)\), \(f = u + v\):
\begin{align*}
    f(x) &= \dfrac{1}{2}\pc{\pa{\alpha + \beta\sqrt{D}}^x + \pa{\alpha - \beta\sqrt{D}}^x} + \dfrac{1}{2\sqrt{D}}\pc{\pa{\alpha + \beta\sqrt{D}}^x - \pa{\alpha - \beta\sqrt{D}}^x} \\
    &= \dfrac{1}{2}\pc{\pa{\alpha + \beta\sqrt{D}}^x + \pa{\alpha + \beta\sqrt{D}}^{-x}} + \dfrac{1}{2\sqrt{D}}\pc{\pa{\alpha + \beta\sqrt{D}}^x - \pa{\alpha + \beta\sqrt{D}}^{-x}} \\
    &= \dfrac{1}{2}\pc{e^{\log\pa{\alpha + \beta\sqrt{D}}x} + e^{-\log\pa{\alpha + \beta\sqrt{D}}x}} + \dfrac{1}{2\sqrt{D}}\pc{e^{\log\pa{\alpha + \beta\sqrt{D}}x} - e^{-\log\pa{\alpha + \beta\sqrt{D}}x}} \\
    &= \cosh\pa{\log\pa{\alpha + \beta \sqrt{D}}x} + \dfrac{1}{\sqrt{D}}\sinh\pa{\log\pa{\alpha + \beta \sqrt{D}}x}.
\end{align*}
For simplicity, let \(y = \log\pa{\alpha + \beta \sqrt{D}}x\). Then
\begin{align*}
    f(x) &= \cosh y + \dfrac{1}{\sqrt{D}}\sinh y \\
    &= \dfrac{2D\cosh y + 2\sqrt{D}\sinh y}{2D} \\
    &= \dfrac{2D\cosh y + \sqrt{4D\pa{\cosh^2 y - 1}}}{2D} \\
    &= \dfrac{2D\cosh y + \sqrt{\pa{2D\cosh y}^2 - 4D\pa{D\cosh^2 y - \cosh^2 y + 1}}}{2D}.
\end{align*}
By applying the quadratic formula, we find that \(f(x)\) satisfies
\begin{align}
    \label{4.1}
    D\pa{f(x)}^2 - \pa{2D\cosh y} f(x) + \pa{D\cosh^2 y - \cosh^2 y + 1} = 0.
\end{align}
Rearranging Equation \eqref{4.1} and solving for \(\cosh y\), we obtain
\[(D-1)\cosh^2 y - \pa{2Df(x)} \cosh y + \pa{D\pa{f(x)}^2 + 1} = 0,\]
which implies
\begin{align}
    \notag
    \cosh y &= \dfrac{2Df(x) \pm \sqrt{\pa{2Df(x)}^2 - 4\pa{D-1}\pa{D\pa{f(x)}^2 + 1}}}{2\pa{D-1}} \\
    \label{4.2}
    &= \dfrac{Df(x) \pm \sqrt{D\pa{f(x)}^2 - D + 1}}{D-1}.
\end{align}
Here the symbol \(\pm\) denotes either the plus \((+)\) or minus \((-)\) sign. Substituting \(x = 0\) into \eqref{4.2}, we find
\[1 = \dfrac{D \pm \sqrt{D - D + 1}}{D-1} = \dfrac{D \pm 1}{D-1},\]
which shows that the minus sign must be chosen. Substitute back \(y = \log\pa{\alpha + \beta \sqrt{D}}x\) into \eqref{4.2} and solve for \(x\), we obtain
\[x = \dfrac{1}{\log\pa{\alpha + \beta \sqrt{D}}}\cosh^{-1}\pa{\dfrac{Df(x) - \sqrt{D\pa{f(x)}^2 - D + 1}}{D-1}}.\]
Finally, replacing \(x\) with \(g(x) := f^{-1}(x)\), we obtain \eqref{2.1}. \\

To obtain \eqref{2.2}, we start from \eqref{2.1} and use the identity \(\cosh^{-1} x = \log \pa{x + \sqrt{x^2 - 1}}\). In particular, for any \(x \geq 1\),
\begin{align*}
    &\cosh^{-1}\pa{\dfrac{Dx - \sqrt{Dx^2 - D + 1}}{D-1}} \\
    &= \log \pa{\dfrac{Dx - \sqrt{Dx^2 - D + 1}}{D-1} + \sqrt{\pa{\dfrac{Dx - \sqrt{Dx^2 - D + 1}}{D-1}}^2 - 1}} \\
    &= \log \pa{\dfrac{Dx - \sqrt{Dx^2 - D + 1} + \sqrt{\pa{Dx - \sqrt{Dx^2 - D + 1}}^2 - \pa{D-1}^2}}{D-1}}
\end{align*}
and
\begin{align*}
    &\sqrt{\pa{Dx - \sqrt{Dx^2 - D + 1}}^2 - \pa{D-1}^2} \\
    &= \sqrt{D^2x^2 + Dx^2 - D + 1 - 2Dx\sqrt{Dx^2 - D + 1} - D^2 + 2D - 1} \\
    &= \sqrt{D}\sqrt{\pa{D+1}x^2 - D + 1 - 2x\sqrt{Dx^2 - D + 1}}.
\end{align*}
\end{proof}
\vspace{10pt}

Expression \eqref{2.2} is rather complicated. From Theorem \ref{thm:2.1}, we observe that the formula for \(Q_{F_D}(\lambda)\) remains valid if we replace the function \(g(\lambda)\) with \(\pe{g\pa{\lambda}}\). This observation motivates us to ask whether there exists a simpler explicit formula for the function \(\pe{g(x)}\) compared to the original expression for \(g(x)\) given in \eqref{2.2}. It turns out that the answer is yes, as we will describe in the following Theorem \ref{thm:2.3}. \\

\begin{thm}
    \label{thm:2.3}
    Let \(g\) be function as described in Theorem \ref{thm:2.2}. Then for any real number \(x \geq 1\),
    \begin{align}
        \label{2.3}
        \pe{g(x)} = \pe{\dfrac{\log\pa{\pe{x}} + C}{\log\pa{\alpha + \beta\sqrt{D}}}},
    \end{align}
    where \(C = \log \pa{\dfrac{2\sqrt{D}}{1 + \sqrt{D}}}\) is a positive real constant. More precisely, there exists a bounded function \(\mu(x): [1,\infty) \to [0,1)\) which vanishes precisely at the points \(f(0), f(1), f(2), \cdots\), and satisfies
    \begin{align}
        \label{2.4}
        g(x) = \pe{\dfrac{\log\pa{\pe{x}} + C}{\log\pa{\alpha + \beta\sqrt{D}}}} + \mu(x)
    \end{align}
    for all \(x \geq 1\).
\end{thm}
\vspace{10pt}

\begin{proof}
We have function \(f: [0,\infty) \to [1,\infty)\) as
\begin{align*}
    f(x) &= \dfrac{1}{2}\pc{\pa{\alpha + \beta\sqrt{D}}^x + \pa{\alpha - \beta\sqrt{D}}^x} + \dfrac{1}{2\sqrt{D}}\pc{\pa{\alpha + \beta\sqrt{D}}^x - \pa{\alpha - \beta\sqrt{D}}^x} \\
    &= \dfrac{1 + \sqrt{D}}{2\sqrt{D}}\pa{\alpha + \beta\sqrt{D}}^x + \dfrac{\sqrt{D} - 1}{2 \sqrt{D}}\pa{\alpha - \beta\sqrt{D}}^x.
\end{align*}
Since \(0 < \dfrac{\sqrt{D} - 1}{2 \sqrt{D}} < 1\) and \(0 < \pa{\alpha - \beta\sqrt{D}}^x = \pa{\alpha + \beta\sqrt{D}}^{-x} \leq 1\), we have
\[0 < \dfrac{\sqrt{D} - 1}{2 \sqrt{D}}\pa{\alpha - \beta\sqrt{D}}^x < 1\]
which implies
\begin{align}
    \notag
    f(x) - 1 &< \dfrac{1 + \sqrt{D}}{2\sqrt{D}}\pa{\alpha + \beta\sqrt{D}}^x < f(x) \\
    \label{5.1}
    \implies \dfrac{2\sqrt{D}}{1 + \sqrt{D}}\pa{f(x) - 1} &< \pa{\alpha + \beta\sqrt{D}}^x < \dfrac{2\sqrt{D}}{1 + \sqrt{D}}\pa{f(x)}.
\end{align}
Focusing on the right-hand side of \eqref{5.1}, we substitute \(x\) with \(\pe{x}\) and then take logarithms, obtaining
\begin{align}
    \label{5.2}
    \pe{x} < \dfrac{\log\pa{f\pa{\pe{x}}} + C}{\log\pa{\alpha + \beta\sqrt{D}}}
\end{align}
where \(C = \log \pa{\dfrac{2\sqrt{D}}{1 + \sqrt{D}}}\). Since \(f\) is strictly increasing, we have \(f\pa{x} \geq f\pa{\pe{x}}\). Because \(f\pa{\pe{x}} \in \Z\), it follows that \(\pe{f\pa{x}} \geq f\pa{\pe{x}}\). Therefore, \eqref{5.2} becomes
\begin{align}
    \label{5.3}
    \pe{x} < \dfrac{\log\pa{\pe{f\pa{x}}} + C}{\log\pa{\alpha + \beta\sqrt{D}}}.
\end{align}
On the other hand, focusing on the left-hand side of \eqref{5.1}, we substitute \(x\) with \(\pe{x} + 1\) and then take logarithms, obtaining
\begin{align}
    \label{5.4}
    \pe{x} + 1 > \dfrac{\log\pa{f\pa{\pe{x} + 1}-1} + C}{\log\pa{\alpha + \beta\sqrt{D}}}
\end{align}
where \(C = \log \pa{\dfrac{2\sqrt{D}}{1 + \sqrt{D}}}\). Since \(f\) is strictly increasing, \(f\pa{\pe{x} + 1} > f\pa{x} \geq \pe{f\pa{x}}\). Since \(f\pa{\pe{x} + 1}\) and \(\pe{f\pa{x}}\) are integers, the strict inequality \(f\pa{\pe{x} + 1} > \pe{f\pa{x}}\) implies \(f\pa{\pe{x} + 1} \geq \pe{f\pa{x}} + 1\). Therefore, \eqref{5.4} becomes
\begin{align}
    \label{5.5}
    \pe{x} + 1 > \dfrac{\log\pa{\pe{f\pa{x}}} + C}{\log\pa{\alpha + \beta\sqrt{D}}}.
\end{align}
By \eqref{5.3} and \eqref{5.5}, we conclude that
\[\pe{x} = \pe{\dfrac{\log\pa{\pe{f\pa{x}}} + C}{\log\pa{\alpha + \beta\sqrt{D}}}}.\]
Finally, replacing \(x\) with \(g(x) := f^{-1}(x)\), we obtain \eqref{2.3}. \\

Since Equation \eqref{2.4} follows from Equation \eqref{2.3}, it remains to show that \(\mu(x) = 0\) if and only if \(x \in \pc{f(0), f(1), f(2), \cdots}\). First suppose that \(\mu(x) = 0\). Then by \eqref{2.4}, we have \(g(x) = f^{-1}(x) \in \Z_{\geq 0}\), which implies that \(x \in \pc{f(0), f(1), f(2), \cdots}\). Conversely, if \(x \in \pc{f(0), f(1), f(2), \cdots}\), then \(g(x) \in \Z_{\geq 0}\). By \eqref{2.4}, it follows that \(\mu(x) = 0\).
\end{proof}
\vspace{10pt}

By combining Theorem \ref{thm:2.1} and Theorem \ref{thm:2.3}, we may conclude the following consequence. \\

\begin{cor}
    Consider the set \(F_D\) as defined in \eqref{1.2}. Then for any real number \(\lambda > 0\),
    \[Q_{F_D}(\lambda) = \left\{
    \begin{aligned}
        &\varnothing \ & \ &\text{if} \ \lambda < 1 \\
        &\bigcup_{i,j \in \pc{1,-1}} \pc{(iu(x), jv(x)) \in \Z^2: x \in \Z_{\geq 0}, x \leq \pe{\dfrac{\log\pa{\pe{\lambda}} + C}{\log\pa{\alpha + \beta\sqrt{D}}}}} \ & \ &\text{if} \ \lambda \geq 1,
    \end{aligned}
    \right.
    \]
    where \(C = \log \pa{\dfrac{2\sqrt{D}}{1 + \sqrt{D}}}\) is a positive real constant. This implies
    \[\left|Q_{F_D}(\lambda)\right| = \left\{
    \begin{aligned}
        &0 \ & \ &\text{if} \ \lambda < 1 \\
        &2 + 4 \pe{\dfrac{\log\pa{\pe{\lambda}} + C}{\log\pa{\alpha + \beta\sqrt{D}}}} \ & \ &\text{if} \ \lambda \geq 1.
    \end{aligned}
    \right.
    \]
\end{cor}
\vspace{10pt}

\section{Further Results}
The Pell equation \eqref{1.1} can be written as \(\pa{x-0}^2 - D\pa{y-0}^2 = 1\), which is the Pell equation centered at the origin. More generally, for any fixed \(a,b \in \Z\), we consider the \textit{shifted Pell equation},
\begin{align}
    \label{6.1}
    \pa{x-a}^2 - D\pa{y-b}^2 = 1.
\end{align}
As an analogue to \eqref{1.2}, we define the set \(F_D^{(a,b)}\) as the set of integer solutions to Equation \eqref{6.1}, i.e.
\begin{align}
    \label{6.2}
    F_D^{(a,b)} &:= \pc{(x,y) \in \Z^2: \pa{x-a}^2 - D\pa{y-b}^2 = 1}.
\end{align}

Then, by Proposition \ref{prop:1.2}, we obtain the following result
\begin{align}
    \label{6.3}
    F_D^{(a,b)} &= \bigcup_{i,j \in \pc{1,-1}} \pc{\pa{iu(x) + a, jv(x) + b} \in \Z^2: x \in \Z_{\geq 0}},
\end{align}
where \(u\) and \(v\) are functions as defined in \eqref{1.3} and \eqref{1.4}, respectively. Similar to the case of integer solutions, the set of real solutions to Equation \eqref{6.1} is given by
\begin{align}
    \label{6.4}
    \bigcup_{i,j \in \pc{1,-1}} \pc{\pa{iu(x) + a, jv(x) + b} \in \R^2: x \in \R, x \geq 0}.
\end{align}
For simplicity, denote each real solution in the set \eqref{6.4} by
\begin{align}
    \label{6.5}
    R_{i,j}(x) := \pa{iu(x) + a, jv(x) + b}.
\end{align}

Notice that the \(x\)-intercepts and \(y\)-intercepts of Equation \eqref{6.1} are \(\pa{a \pm \sqrt{1+Db^2},0}\) and \(\pa{0, b \pm \sqrt{\dfrac{a^2 - 1}{D}}}\), respectively (there is no \(y\)-intercept if \(a = 0\)). \\

For any \(k \in \{1,2,3,4\}\), let \(\mathcal{H}_k\) denote the set of all real pairs in the \(k\)-th quadrant of the \(\R^2\) plane, including the boundary axes. The following proposition describes the location of the real solutions \(R_{i,j}(x)\) in the \(\R^2\) plane for all sufficiently large \(x \geq 0\): \\

\begin{prop}
    \label{prop:6.1}
    Let \(R_{i,j}(x)\) be a real solution as defined in \eqref{6.5}. Then \(R_{1,1}(x) \in \mathcal{H}_1\), \(R_{-1,1}(x) \in \mathcal{H}_2\), \(R_{1,-1}(x) \in \mathcal{H}_3\) and \(R_{-1,-1}(x) \in \mathcal{H}_4\) for all sufficiently large \(x \geq 0\).
\end{prop}
\vspace{10pt}

\begin{proof}
    Without loss of generality, let \(a\) and \(b\) be non-negative integers (the cases \(a < 0\) or \(b < 0\) can be shown similarly). Then \(R_{1,1}(x) \in \mathcal{H}_1\) for all \(x \in [0,\infty)\). For the case of \(a > 0\), there exist \(x\)-intercepts and \(y\)-intercepts to the equation \(\pa{x-a}^2 - D\pa{y-b}^2 = 1\). In particular, there exist \(p, q, r \in [0,\infty)\) such that 
    \begin{enumerate}
        \item \(R_{-1,1}(p) = \pa{0, b + \sqrt{\dfrac{a^2 - 1}{D}}} \in \mathcal{H}_2\),
        \item \(R_{-1,-1}(q) = \left\{
        \begin{aligned}
            &\pa{0, b - \sqrt{\dfrac{a^2 - 1}{D}}} \ & \ &\text{if} \ \pa{0, b - \sqrt{\dfrac{a^2 - 1}{D}}} \in \mathcal{H}_3 \\
            &\pa{a - \sqrt{1+Db^2},0} \ & \ &\text{if} \ \pa{a - \sqrt{1+Db^2},0}\in \mathcal{H}_3 \hspace{10pt} \text{and} \\
        \end{aligned}
        \right. \)
        \item \(R_{1,-1}(r) = \pa{a + \sqrt{1+Db^2},0} \in \mathcal{H}_4\). \\
    \end{enumerate}
    Since functions \(u\) and \(v\) are strictly increasing, this implies \(R_{-1,1}(x) \in \mathcal{H}_2\) for all \(x \in [p, \infty)\), \(R_{-1,-1}(x) \in \mathcal{H}_3\) for all \(x \in [q, \infty)\) and \(R_{1,-1}(x) \in \mathcal{H}_4\) for all \(x \in [r, \infty)\). \\

    For the case of \(a = 0\), we have \(R_{-1,1}(x) \in \mathcal{H}_2\) for all \(x \in [0,\infty)\). There exist \(x\)-intercepts (but no \(y\)-intercept) to Equation \eqref{6.1}. By similar argument, there exist \(q, r \in [0,\infty)\) such that \(R_{-1,-1}(x) \in \mathcal{H}_3\) for all \(x \in [q, \infty)\) and \(R_{1,-1}(x) \in \mathcal{H}_4\) for all \(x \in [r, \infty)\).
\end{proof}
\vspace{10pt}

The following proposition helps establish the setting of Theorem \ref{thm:6.3}: \\
\begin{prop}
    \label{prop:6.2}
    Let \(a \neq 0\). Then
    \[|a| + \sqrt{1+Db^2} = \max\pc{\pd{a \pm \sqrt{1+Db^2}}, \pd{b \pm \sqrt{\dfrac{a^2 - 1}{D}}}}.\]
\end{prop}
\vspace{10pt}

\begin{proof}
    Notice that
    \[|a| + \sqrt{1+Db^2} = \max\pc{\pd{a \pm \sqrt{1+Db^2}}} \hspace{10pt} \text{and} \hspace{10pt} |b| + \sqrt{\dfrac{a^2 - 1}{D}} = \max\pc{\pd{b \pm \sqrt{\dfrac{a^2 - 1}{D}}}},\]
    so it suffices to show that \(|a| + \sqrt{1+Db^2} \geq |b| + \sqrt{\dfrac{a^2 - 1}{D}}\). This is equivalent to \(\sqrt{a^2} + \sqrt{1+Db^2} \geq \sqrt{b^2} + \sqrt{\dfrac{a^2 - 1}{D}}\), which must be true since \(\sqrt{a^2} > \sqrt{\dfrac{a^2 - 1}{D}}\) and \(\sqrt{1+Db^2} > \sqrt{b^2}\).
\end{proof}
\hspace{10pt}

We now explicitly enumerate the set \(Q_{F_D^{(a,b)}}(\lambda)\) for any sufficiently large real number \(\lambda > 0\), as follows: \\
\begin{thm}
    \label{thm:6.3}
    Consider the set \(F_D^{(a,b)}\) as defined in \eqref{6.2}. Set
    \[K := |a| + \max\pc{|b| + 1, \sqrt{1 + Db^2}}.\]
    Then for any real number \(\lambda \geq K\),
    \[Q_{F_D^{(a,b)}}(\lambda) = \bigcup_{i,j \in \pc{1,-1}} \pc{\pa{iu(x) + a, jv(x) + b} \in \Z^2: x \in \Z_{\geq 0}, x \leq g\pa{\lambda - ia - jb}}, \]
    where \(g\) is the function as described in Theorem \ref{thm:2.3}. This implies for any real number \(\lambda \geq K\),
    \[\left| Q_{F_D^{(a,b)}}(\lambda) \right| = 2 + \sum_{i,j \in \{1,-1\}} \pe{g\pa{\lambda - ia - jb}}.\]
\end{thm}
\vspace{10pt}

\begin{proof}
By referring to the expression \(F_D^{\pa{a,b}}\) in \eqref{6.3}, we obtain
\[Q_{F_D^{(a,b)}}(\lambda) = \bigcup_{i,j \in \pc{1,-1}} \pc{\pa{iu(x) + a, jv(x) + b} \in \Z^2: x \in \Z_{\geq 0}, \pd{iu(x) + a} + \pd{jv(x) + b} \leq \lambda}.\]

For simplicity, let \(\overline{f}_{i,j} : [0, \infty) \to [0,\infty)\) with \(\overline{f}_{i,j}(x) := \pd{iu(x) + a} + \pd{jv(x) + b}\). Then, it remains to show: \\
\begin{enumerate}[label=(\Alph*)]
    \item \label{A} \(\overline{f}_{i,j}(x) \leq \lambda\) for all non-negative integers \(x \leq g(\lambda - ia - jb)\) and \\
    \item \label{B} \(\overline{f}_{i,j}(x) > \lambda\) for all integers \(x > g(\lambda - ia - jb)\). \\
\end{enumerate}

\begin{rem}
    \label{rem:6.4}
    Since \(\overline{f}_{i,j}\) is the sum of two absolute functions, and since \(u\) and \(v\) are both strictly increasing, there exist two real numbers \(0 \leq x_0 \leq x_1 < \infty\) such that the function \(\overline{f}_{i,j}\) is monotone (not necessary strictly) decreasing on the interval \([0,x_0)\), monotone (not necessary strictly) increasing on \([x_0,x_1)\) and strictly increasing on \([x_1, \infty)\). In other words, if \(\overline{f}_{i,j}(x) < \overline{f}_{i,j}(y)\) for some \(0 \leq x < y < \infty\), then the function \(\overline{f}_{i,j}\) must be monotone increasing for all points after \(y\).
\end{rem}
\vspace{10pt}

Before proving \ref{A} and \ref{B}, we require a few more lemmas. \\

\begin{lem}
    \label{lem:6.5}
    For any \(i,j \in \{1,-1\}\), \(\overline{f}_{i,j}(0) \leq K\).
\end{lem}
\vspace{10pt}

\begin{proof}
By the triangle inequality, \(\overline{f}_{i,j}(0) = |a \pm 1| + |b| \leq |a| + |b| + 1 \leq K\).
\end{proof}
\vspace{10pt}

Since \(\overline{f}_{i,j}\) is eventually strictly increasing, by Lemma \ref{lem:6.5}, there exists some \(x \in [0,\infty)\) such that \(\overline{f}_{i,j}(x) = K\). We then define
\begin{align}
    \label{6.6}
    x_1 := \sup\pc{x \in [0, \infty): \overline{f}_{i,j}(x) = K}.
\end{align}
where \(x_1\) depends on \(i,j \in \{1,-1\}\). \\

\begin{lem}
    \label{lem:6.6}
     Let \(x_1\) be the real value as defined in \eqref{6.6}. Then for any \(x \in [0,x_1]\), \(\overline{f}_{i,j}(x) \leq K\).
\end{lem}
\vspace{5pt}

\begin{proof}
    Notice that the statement holds for \(x = 0\) (from Lemma \ref{lem:6.5}) and for \(x = x_1\) (from \eqref{6.6}). Assume that there exists \(y \in (0,x_1)\) such that \(\overline{f}_{i,j}(y) > K\). By Lemma \ref{lem:6.5}, this implies \(\overline{f}_{i,j}(0) < \overline{f}_{i,j}(y)\). Hence, by Remark \ref{rem:6.4}, the function \(\overline{f}_{i,j}\) must be monotone increasing for all points after \(y\). However, at the point \(x_1\) we have \(\overline{f}_{i,j}(y) > K = \overline{f}_{i,j}(x_1)\) from \eqref{6.6}, a contradiction.
\end{proof}
\vspace{10pt}

On the other hand, lemma below would implies \(\overline{f}_{i,j}(x) \geq K\) for all \(x \in [x_1, \infty)\): \\
\begin{lem}
    \label{lem:6.7}
    Let \(x_1\) be the real value as defined in \eqref{6.6}. Then for any \(i,j \in \{1,-1\}\), the function \(\overline{f}_{i,j} : [x_1, \infty) \to [K,\infty)\) is bijective and strictly increasing. In particular, we have the function \(\overline{f}_{i,j} : [x_1, \infty) \to [K,\infty)\) as
    \begin{align}
        \label{6.7}
        \overline{f}_{i,j}(x) = f(x) + ia + jb,
    \end{align}
    where \(f = u + v\).
\end{lem}
\vspace{10pt}

\begin{proof}
     Without loss of generality, let \(a\) and \(b\) be non-negative integers (the cases \(a < 0\) or \(b < 0\) can be shown similarly). By Proposition \ref{prop:6.1}, there exists \(p,q,r \in [0,\infty)\) (if \(a = 0\) then \(p=0\)) such that \(R_{1,1}(x) \in \mathcal{H}_1\) for all \(x \in [0, \infty)\), \(R_{-1,1}(x) \in \mathcal{H}_2\) for all \(x \in [p, \infty)\), \(R_{-1,-1}(x) \in \mathcal{H}_3\) for all \(x \in [q, \infty)\) and \(R_{1,-1}(x) \in \mathcal{H}_4\) for all \(x \in [r, \infty)\). \\
    
    Notice that for any \(i,j\in \{1,-1\}\), we have \(\|R_{i,j}(x)\| = \overline{f}_{i,j}(x)\) for all \(x \in [0,\infty)\), where \(\|\cdot\|_1\) denotes the \(\mathcal{L}^1\)-norm. Since \(R_{1,1}(x) \in \mathcal{H}_1\) for all \(x \in [0, \infty)\), it follows that
    \[\|R_{1,1}(x)\| = \pa{u(x) + a} + \pa{v(x) + b},\]
    which implies \(\overline{f}_{1,1}(x) = f(x) + a + b\) for all \(x \in [0,\infty)\). By similar argument, we have \(\overline{f}_{-1,1}(x) = f(x) - a + b\) for all \([p,\infty)\), \(\overline{f}_{-1,-1}(x) = f(x) - a - b\) for all \([q,\infty)\) and \(\overline{f}_{1,-1}(x) = f(x) + a - b\) for all \([r,\infty)\). Therefore, we have functions
    \begin{align*}
        \overline{f}_{1,1}:[0,\infty) &\to \Bigl[\overline{f}_{1,1}(0), \infty \Bigr) \\
        \overline{f}_{-1,1}:[p,\infty) &\to \Bigl[\overline{f}_{-1,1}(p), \infty \Bigr) \\
        \overline{f}_{-1,-1}:[q,\infty) &\to \Bigl[\overline{f}_{-1,-1}(q), \infty \Bigr) \\
        \overline{f}_{1,-1}:[r,\infty) &\to \Bigl[\overline{f}_{1,-1}(r), \infty \Bigr)
    \end{align*}
    with \(\overline{f}_{i,j}(x) = f(x) + ia + jb\) for any \(i,j \in \{1,-1\}\), where each \(\overline{f}_{i,j}\) is bijective and strictly increasing. \\

    It remains to show that \(p,q,r \leq x_1\). We will only show for \(p \leq x_1\), as \(q \leq x_1\) and \(r \leq x_1\) can be shown using a similar approach. If \(a = 0\), then \(p = 0 \leq x_1\), so without loss of generality, let \(a > 0\). Assume that \(p > x_1\). Since \(R_{-1,1}(p)\) is the \(y\)-intercept (as described in the proof of Proposition \ref{prop:6.1}), Proposition \ref{prop:6.2} implies that \(\|R_{-1,1}(p)\| \leq K\). Equivalently, \(\overline{f}_{-1,1}(p) \leq K = \overline{f}_{-1,1}(x_1)\). If \(\overline{f}_{-1,1}(p) = K = \overline{f}_{-1,1}(x_1)\), then it contradicts the definition of \(x_1\) in \eqref{6.6}, so we must have \(\overline{f}_{-1,1}(p) < K = \overline{f}_{-1,1}(x_1)\). But \(\overline{f}_{-1,1}\) is bijective and strictly increasing on the domain \([p,\infty)\), so there exists \(\overline{x} > p > x_1\) such that \(\overline{f}_{-1,1}\pa{\overline{x}} = K\), which contradicts the definition of \(x_1\) in \eqref{6.6}.
\end{proof}
\vspace{10pt}

We are ready to show \ref{A} and \ref{B}: \\

\textit{Proof of \emph{\ref{A}}.} Suppose that there exists a non-negative integer \(x \leq g(\lambda - ia - jb)\) such that \(\overline{f}_{i,j}(x) > \lambda \geq K\). By the contrapositive of Lemma \ref{lem:6.6}, we must have \(x \in (x_1, \infty)\). Since \(f\) is strictly increasing and \(x \leq g(\lambda - ia - jb)\), it follows that \(f(x) + ia + jb \leq \lambda\). By \eqref{6.7}, this is equivalent to \(\overline{f}_{i,j}(x) \leq \lambda\), a contradiction. \\

\textit{Proof of \emph{\ref{B}}.} Pick any integer \(x > g(\lambda - ia - jb)\). Since \(f\) is strictly increasing, we have \(f(x) + ia + jb > \lambda \geq K\). By Lemma \ref{lem:6.7}, \(\overline{f}_{i,j}\) is surjective, so there exists \(y \in [x_1, \infty)\) such that \(\overline{f}_{i,j}(y) = f(x) + ia + jb\). Because \(\overline{f}_{i,j}\) is injective, it follows from \eqref{6.7} that \(y = x\) and hence \(\overline{f}_{i,j}(x) > \lambda\).  \\

For the number of elements, \(\left| Q_{F_D^{(a,b)}}(\lambda) \right|\), we have
\begin{align*}
    Q_{F_D^{(a,b)}}(\lambda) &= \bigcup_{i,j \in \pc{1,-1}} \pc{\pa{iu(x) + a, jv(x) + b} \in \Z^2: x \in \Z_{\geq 0}, x \leq g\pa{\lambda - ia - jb}} \\
    &= \pc{(a+1,b), (a-1,b)} \\
    &\hspace{30pt} \cup \bigcup_{i,j \in \pc{1,-1}} \pc{\pa{iu(x) + a, jv(x) + b} \in \Z^2: x \in \N, x \leq \pe{g\pa{\lambda - ia - jb}}}
\end{align*}
for any \(\lambda \geq K\), which implies \(\left| Q_{F_D^{(a,b)}}(\lambda) \right| = 2 + \displaystyle\sum_{i,j \in \{1,-1\}} \pe{g\pa{\lambda - ia - jb}}\).
\end{proof}

\section{Conclusion and Future Work}
We have given an explicit enumeration of integer solutions to the Pell equation \eqref{1.1} within the square \(|x| + |y| \leq \lambda\) for any \(\lambda > 0\) and extended our results to the shifted Pell equation \eqref{6.1} for all sufficiently large \(\lambda > 0\). Unlike previous studies that focus on the asymptotic behavior of fundamental units as \(D\) varies, our results provide exact counts for a fixed \(D\), thereby complementing existing literature. \\

Future directions include extending the method to other bounded regions such as discs, rectangles, or convex bodies and applying similar ideas to more general quadratic Diophantine equations such as \(Ax^2 + Bxy + Cy^2 = 1\) or to higher-dimensional analogues. In addition, one may consider the broader family of Pell-type equations \(x^2 - Dy^2 = n\) with \(n \in \Z\) and study the distribution of their integer solutions within bounded regions. This setting raises new challenges, since the solvability depends on congruence conditions for \(n\) and the solution structure may vary significantly with \(n\). \\

These directions suggest that the study of Pell-type equations within bounded regions is far from complete. Our work provides an initial step in this line of research and we anticipate that further developments will yield deeper insights into the fine structure of Diophantine equations.

\vskip 20pt\noindent {\bf Acknowledgement.}
The research in this paper was conducted during the first author’s master’s studies at Universiti Kebangsaan Malaysia (UKM). \\

\end{document}